\documentclass[12pt]{amsart}
\usepackage{amsmath,amsfonts,mathrsfs,enumerate,graphics,color,amssymb,txfonts}
\usepackage[all]{xy}

\theoremstyle{theorem}
\newtheorem{theorem}{Theorem}[section]

\newtheorem{corollary}[theorem]{Corollary}

\newtheorem{prop}[theorem]{Proposition}
\newtheorem{lem}[theorem]{Lemma}
\theoremstyle{rem}

\theoremstyle{definition}

\newtheorem{rema}[theorem]{Remark}

\numberwithin{figure}{section}
\numberwithin{equation}{section}

\newcommand{\f}{\mathfrak f}

\newcommand{\ifff}{if and only  if }

\newcommand{\ZZ}{\mathbb Z}
\newcommand{\PP}{\mathbb P}

\newcommand{\Ii}{\mathscr I}

\newcommand{\Qq}{\mathcal Q}

\DeclareMathOperator{\coker}{coker}
\DeclareMathOperator{\Ext}{Ext}

\newcommand{\Py}{P_y}\newcommand{\Px}{P_x}
\newcommand{\ra}{\longrightarrow}

\newcommand{\C}{\mathcal C}

\newcommand{\Hh}{\mathcal H}

\newcommand{\cI}{{\mathcal I}}

\DeclareMathOperator{\Ann}{Ann}

\begin{document}
\title[The equations of space curves on a quadric]{The equations of
  space curves on a quadric}%
\author[Roberta Di Gennaro]{Roberta Di Gennaro$^1$}
\author{Uwe Nagel}
\thanks{$^1$ Supported by Universit\`{a} degli Studi di Napoli \lq\lq
  Federico II\rq\rq\ in the framework of the \lq\lq Programma di
  Scambi
 Internazionali per la Mobilit\`{a} di breve durata di docenti,
 studiosi e ricercatori\rq\rq}
\address{Dipartimento di Matematica e Applicazioni \lq\lq R. Caccioppoli\rq
\rq, Universit\`{a} degli Studi di Napoli \lq\lq Federico
II\rq\rq, Complesso Universitario Monte S. Angelo - 80126 Napoli,
Italy}
\email{digennar@unina.it}
\address{Department of Mathematics,
 University of Kentucky,
 715 Patterson Office Tower,
 Lexington, KY 40506-0027,
 USA
}
\email{uwenagel@ms.uky.edu}
\par \vspace{.5cm}
\date{}


\begin{abstract}
The homogenous ideals of curves in a double plane have been studied by
Chiarli, Greco,
Nagel. Completing this work we describe the equations of any curve
that is contained in
some quadric. As a consequence, we classify the Hartshorne-Rao modules
of such curves.
\end{abstract}

\maketitle

\section{Introduction}
%

The goal of this note is to study the equations of a curve $C \subset
\PP^3 = \PP^3_K$ that is contained in some quadric $\Qq$. By a curve, we mean a
pure one-dimensional locally Cohen-Macaulay subscheme (i.e.\ without
zero-dimensional components.)  We assume that the field $K$ is
algebraically closed.

If $C \subset {\mathcal Q}$ is arithmetically Cohen-Macaulay, then, by Dubreil's Theorem,
$C$ is defined by at most 3 equations. The converse is also true by a well-known result
of Evans and Griffith (\cite[Theorem 2.1]{EG}).  Denoting by $\mu (I_C)$ the number of
minimal generators of the homogeneous ideal of $C$, this gives:

\begin{prop} \label{prop-aCM}
If $C \subset \PP^3$ is a curve lying on some quadric, then $C$ is
arithmetically
Cohen-Macaulay if and only if $\mu (I_C) \leq 3$.
\end{prop}

In this case, $C$ is either a complete intersection or its ideal $I_C$ is generated by
the $2$-minors of a $2 \times 3$ matrix.
\par\vspace{3mm}
In order to discuss the ideal of $C \subset \Qq$ when $C$ is not arithmetically
Cohen-Macaulay, we take the rank of the quadric $\Qq$ into account.
If $\Qq$ has rank one, then $\Qq$ is not reduced, i.e.\ $\Qq$ is a
double plane $2 H$. Building on the work of Hartshorne and Schlesinger
\cite{HS}, the homogeneous ideal of $C$ has been described by Chiarli,
Greco, and Nagel in \cite{CGN_doubleplane}.

In case $\mathcal Q$  has rank 2, the curves on $\mathcal Q$ have been
studied by
Hartshorne in \cite{H-gen-div}, Section 5,  from the point of view of
generalized divisors. The quadric $\Qq$ is a union of two distinct
hyperplanes. In Section \ref{sec-redquadric} we establish the following
characterization of the equations of curves on such a reducible quadric:

\begin{theorem}
\label{thm-red-quad}
 Let $C$ be a curve on a reducible reduced quadric $\mathcal Q=\Hh \cup
\Hh'$. Then $C$ is not arithmetically Cohen-Macaulay \ifff $\mu(I_C)
=4$. In this case
$C$ is (up to a change of coordinates) defined by
$$
I_C=\left(xy, x^2A+ xhF, y^2B + yhG, xAF+yBG+hFG\right)
$$
where $x,y,F,G$ is a regular sequence, $h \in K[z, t]$, and  $y \nmid A$, $x \nmid B$ if
$h = 0$.
\end{theorem}

Using this result, we find the minimal free resolution of such a curve
and, finally, determine its Hartshorne-Rao module $M_\C :=
\oplus_{j\in \ZZ} H^1(\Ii_C(j))$  in Corollary \ref{cor-ModRao-red}.
\par\vspace{3mm}
If $\mathcal Q$ has rank 3, then it is cone. Thus all curves on it are arithmetically
Cohen-Macaulay (cf.\ \cite[Example 5.2]{fe}). This leaves us with the case when $\Qq$ has
rank 4, i.e.\ $\Qq$ is a smooth quadric.

Curves on a smooth quadric are investigated in Section
\ref{sec-smoothquadric}. We use results from liaison theory (cf.\
\cite{migliore} for a comprehensive introduction).
The Lazarsfeld-Rao property of space curves
says that each curve $C$ can be obtained from a so-called minimal
curve in its even liaison class by a sequence of basic double links
and, possibly, a flat deformation. The deformation can be avoided if
we replace basic double links by ascending elementary biliaisons
(\cite{strano}). Thus, it suffices to consider minimal curves. It is
easy to see that any such curve $C$ is linearly equivalent to $d \f$
for some $d \geq 1$ where $\f$ is a line on $\Qq$. This allows us to
compute  explicitly the Hartshorne-Rao module of $C$ in Corollary
\ref{cor-hr-mod}. It improves the results in  \cite{GM-02}. Then we
determine the equations of $C$:

\begin{theorem}\label{thm-dec-smooth}
Let $C\subset \Qq$ be a curve in $|d\f|$. Then there are distinct
lines $L_1,\ldots,L_m$
in the same ruling of $\Qq$ as $\f$ such that $C$ is
 defined by:
$$
I_C=I_{\Qq} + I_{L_1}^{d_1} \cdot  \ldots \cdot  I_{L_m}^{d_m}
$$
where  $d_i\geqslant 1$ and $d_1 + \ldots + d_m = d$.
\end{theorem}

In the final Section \ref{sec-HR} we summarize the results about the
Hartshorne-Rao modules of curves on a quadric. Each such module
determines an even liaison class. In each class we exhibit a minimal curve.

\section{Curves on the union of two distinct planes}\label{sec-redquadric}
In this section, we consider curves on a reducible quadric $\Qq$,
i.e.\ $\Qq$ is the
union of two distinct planes in $\PP^3$. We denote by $R :=
K[x,y,z,t]$ the coordinate ring of $\PP^3$.
Without loss of generality we
may assume that
$\Qq$ is defined by $q = xy$. We will prove Theorem \ref{thm-red-quad}
and determine the
Hartshorne-Rao module of the curves.

We first show that the equations of the curves have the shape as
predicted by Theorem
\ref{thm-red-quad}.

\begin{lem} \label{lem-red-quad-necc}
Let $C$ be a curve on $\{x y = 0\}$, then
\begin{itemize}
\item[(a)] $\mu (I_C) \leq 4$. \\
\item[(b)] $C$ is not arithmetically Cohen-Macaulay \ifff $\mu(I_C) =4$.
In this case $C$ is defined by the homogeneous ideal
$$
I_C=\left(xy, x^2A+ xhF, y^2B + yhG, xAG+yBF+hFG\right)
$$
where $h, F, G \in K[z,t]$, $F, G \notin K$, $A \in K[x,z,t]$, $B \in K[y,z,t]$, and $A B
\neq 0$ in case $h =0$.
\end{itemize}
\end{lem}
\begin{proof}   As in \cite{CGNeghs} we will utilize residual sequences. We will use the
following notation: $R_x=K[y,z,t], R_y=K[x,z,t]$. By Proposition \ref{prop-aCM}, we may
assume that $C$ is not arithmetically Cohen-Macaulay.

Let $C_x \subset \Hh := \{x = 0\}$ and $C_y \subset \Hh' := \{y = 0\}$ be the planar
curves defined by
$$I_{C_y}=I_{C}:x=(y,\Py)$$ and
$$I_{C_x}=I_{C}:y=(x,\Px),$$ with homogeneous polynomials $P_x \in R_x \setminus K$ and $P_y \in
R_y \setminus K$. Denote by $D_x$ the one-dimensional part of $C \cap \Hh$ and by $Z_x$
the residual subscheme to $P_x$ in $C \cap \Hh$. $Z_x$ is not empty because $C$ is not
arithmetically Cohen-Macaulay. Since $D_x$ is the largest planar subcurve of $C \cap
\Hh$, we get $C_x \subset D_x$. Thus, $D_x$ is defined by an ideal $(x, f P_x)$ for some
homogeneous polynomial $0 \neq f \in R_x$. Moreover, \cite[Lemma 2.8]{CGNeghs}  provides
that
$$
Z_x \subset \Hh \cap C_y.
$$
Thus $Z_x$ is defined by an ideal $(x, y, Q_x)$ for some $ Q_x \in K[z,t] \setminus K$.
Hence, the residual sequence of $C$ with respect to $\Hh$ reads as:
\begin{equation} \label{eq-res-seq}
0 \ra (y,P_y)(-1) \stackrel{x}{\ra} I_{\C} \ra P_x\,f\,(y,Q_x) \cdot R_x \ra 0.
\end{equation}
Similarly, we get for the residual sequence with respect to $\Hh'$:
$$
0 \ra (x,P_x)(-1) \stackrel{y}{\longrightarrow} I_{\C} \ra P_y\,k\,(x,Q_y) \cdot R_y \ra
0
$$
where  $0 \neq k \in R_y$ and $Q_y \in K[z,t] \setminus K$.
These sequences imply that we can write the ideal of $C$ as:
\begin{align*}
 I_{\C} & = (xy,\ xP_y,\ yP_x\,f + xA_1 , \ P_x\,f\, Q_x +xA_2) &
      \text{for some } A_1,A_2 \in R_y \\
     & = (xy,\ yP_x,\ xP_y\,k+ y B_1 ,\ P_y \,k\, Q_y +yB_2)&
     \hspace{.5cm} \text{for some } B_1,B_2 \in R_x
\end{align*}
This shows in particular that $\mu (I_C) \leq 4$. Furthermore, it follows that $y B_1 \in
I_C$, thus $B_1$ is in $I_C : y = (x, P_x)$. Since $B_1$ and $P_x$ are in $R_x$, we see
that $P_x$ must divide $B_1$. Analogously, we get that $P_y$ divides $A_1$. Hence, we can
rewrite the ideal of $C$ as
\begin{eqnarray*}
I_{\C} & = & (xy,\ xP_y,\ yP_x, \ P_x\,f\, Q_x +xA_2) \\
& = & (xy,\ xP_y,\ yP_x,\ P_y \,k\, Q_y +yB_2).
\end{eqnarray*}
Comparing with the residual sequence \eqref{eq-res-seq} we obtain:
\begin{eqnarray*}\label{mod x}
I_{C}+x R & = & x R+ (yP_x, P_x\,f\, Q_x)\\
        & = &xR+P_x f\, (y, Q_x ).
\end{eqnarray*}
Since the degree of $Q_x$ is at least one, $f$ must be a constant and we may assume that
$f = 1$. Analogously, we get $k = 1$ without loss of generality, thus
\begin{equation} \label{eq-ideal}
I_C = (xy,\ xP_y,\ yP_x,\ P_y\, Q_y +yB_2)
\end{equation}
We now distinguish two cases.

\noindent \textsc{Case 1.} Assume that the curves $C_x$ and $C_y$ have a common
component. This component must be the line $\Hh \cap \Hh'$. If follows that there are
polynomials $0 \neq A \in R_y$ and $0 \neq B \in R_x$ such that
$$
P_x = y B \quad {\rm and} \quad P_y = x A.
$$
Thus, the ideal of $C$ reads as
$$
I_C = (xy,\ x^2 A,\ y^2 B,\ x\,A\, Q_y + y\,B_2).
$$
Another comparison with the residual sequence \eqref{eq-res-seq} provides
\begin{eqnarray*}\label{mod x}
I_{C}+x R & = & x R+ (y^2 B, y B_2)\\
        & = & x R + y B\, (y, Q_x ).
\end{eqnarray*}
Since $B, B_2, Q_x$ are in $R_x$, we conclude that $B$ divides $B_2$, i.e.\ $B_2 = B F$
for some $ F \in R_x \setminus K$. Setting $G := Q_y$, we get
$$
I_C = (xy,\ x^2A,\ y^2B,\ xAG+yBF).
$$
It follows that we may assume $F, G \in K[z,t] \setminus K$. Thus, $I_C$ is of the
required form
(with $h = 0$). \\[3pt]
\textsc{Case 2.} Assume that $C_x$ and $C_y$ do not have a common component. Then we get
that $C \subset C_x \cup C_y$ and $\deg C = \deg C_x + \deg C_y = \deg (C_x \cup C_y)$.
Since all these curves are of pure dimension one, we conclude that $C = C_x \cup C_y$,
i.e.\
$$
I_C = (x, P_x) \cap (y, P_y)
$$
where $\dim R/(x, y, P_x, P_y) \leq 1$. \\[3pt]
\textsc{Case 2.1.} Assume $\dim R/(x, y, P_x, P_y) = 0$. Then $C$ is the disjoint union
of the planar curves $C_x$ and $C_y$. Write $P_x \in R_x$ as $P_x = y B + G$  for some $G
\in K[z, t] \setminus K$ and some $B \in R_x$ and, similarly, $P_y = x A + F$ for some $F
\in K[z, t] \setminus K$ and some $A \in R_y$. Then we get using also Identity
\eqref{eq-ideal}:
\begin{eqnarray*}
I_C & = &(x, y B + G) \cap (y, x A + F) \\
& = & (x y,\ x^2 A + x F,\ y^2 B + y G,  x A G + y B F + F G)
\end{eqnarray*}
showing that $I_C$ has the required form (with $h = 1$).
\\[3pt]
\textsc{Case 2.2.} Assume that $\dim R/(x, y, P_x, P_y) = 1.$ Then the saturation of $(x,
y, P_x, P_y)$ is of the form $(x, y, h)$ for some polynomial $h \in K[z, t] \setminus K$.
It follows that
$$
P_x = y B + h G \quad {\rm and} \quad P_y = x A + h F
$$
for some regular sequence $F, G \in K[z, t]$ and some $B \in R_x$, $A \in R_y$. Thus, we
obtain:
$$
(xy, x^2A+ xhF, y^2B + yhG, xAG+yBF+hFG) \subset (x, y B + h G) \cap (y, x A + h F) =
I_C.
$$
Comparing with the identity \eqref{eq-ideal} we conclude that the above ideals are equal.
Thus,
 the proof is complete.
\end{proof}

The above lemma specifies necessary conditions on the homogeneous ideals of curves on a
reducible quadric $\Qq$. In order to find sufficient conditions we determine the minimal
free resolution of the ideals. This will also allow us to compute the Hartshorne-Rao
module of the curves.

\begin{lem} \label{lem-true-curve} Consider the following homogeneous ideal in R:
$$
J=\left(xy, x^2A+ xhF, y^2B + yhG, xAG+yBF+hFG\right)
$$
where $h, F, G \in K[z,t]$, $F, G \notin K$, $A \in K[x,z,t]$, $B \in K[y,z,t]$,  and $A
B \neq 0$ in case $h =0$. Denote by $d_F, d_G, d_h$ the degree of $F, G$, and $h$,
respectively. Then $J$ defines a 1-dimensional subscheme $C \subset \PP^3$ of degree $d =
2 d_h + d_F + d_G$.

Moreover, $J$ has the following minimal free graded resolution:
\begin{equation}\label{res_ideal_red}
\xymatrix{
0 \ar[r] & F_3 \ar[r]_-{\phi_3}& F_2
 \ar[r]_-{\phi_2
} &F_1
 \ar[r]_-{\phi_1}&
 J
 \ar[r] & 0
 }
\end{equation}
where
\begin{eqnarray*}
F_1 & = &R(-2) \oplus R(-d_F-d_h-1) \oplus R(-d_G -d_h-1) \oplus R(-d_F-d_G-d_h) \\
F_2 & = &R(-d_F-d_h-2) \oplus R(-d_G -d_h-2) \oplus R(-d_F-d_G-d_h-1)^2 \\
F_3 & = & R(-d_F-d_G-d_h-2)
\end{eqnarray*}
and, by identifying the maps with its matrices,
\begin{eqnarray*}
\phi_1 & = &\begin{bmatrix}
xy, x^2A+ xhF, y^2B + yhG, xAG+yBF+hFG\end{bmatrix}, \\[1ex]
  \phi_2& = &
\begin{bmatrix}
xA + hF & yB + hG & AG & BF \\
  -y& 0&0&G \\
  0&-x&F&0\\
  0&0&-y&-x
\end{bmatrix}, \quad
\phi_3 \ = \ \begin{bmatrix}
 F\\-G\\-x\\y
\end{bmatrix}.
\end{eqnarray*}
Hence, $C$ is a curve if and only if $x,y,F,G$ is a regular sequence.
\end{lem}
\begin{proof}
It is immediate to verify that the Sequence \eqref{res_ideal_red} is a complex. In order
to check its exactness we use the Buchsbaum-Eisenbud criterion (\cite[Theorem 20.9]{E}).
Since $F, G \neq 0$ are in $K[z, t]$, the entries of $\phi_3$ generate an ideal whose
codimension is at least three. It remains to show that the ideal $I_3(\phi_2)$ generated
by the 3-minors of $\phi_2$ contains a regular sequence of length two. Clearly, we have
$x^2 y \in I_3(\phi_2)$. Moreover, the determinant of the matrix obtained from $\phi_2$
by deleting row 3 and column 1 is $y G ( y B + h G) \in I_3(\phi_2)$. It is not divisible
by $x$. Similarly, deleting row 2 and column 3 of $\phi_2$ we get $x^2 (xA + h F) \in
I_3(\phi_2)$ which is not divisible by $y$. It follows that $(x^2 y,\ y G ( y B + h G),\
x^2 (xA + h F)) \subset I_3(\phi_2)$ has codimension two.

Since the non-zero entries of the matrices $\phi_1, \phi_2, \phi_3$ have positive degree,
the Resolution \eqref{res_ideal_red} is minimal.

Finally, $C$ is curve if and only if $\Ext^2_R (I_C, R)$ has finite length. But the
resolution of $I_C$ provides
\begin{equation} \label{eq-ext}
\Ext^2_R (I_C, R) \cong \left ({R}/{(x,y,F,G)}\right ) (d_F + d_G + d_h + 2).
\end{equation}
Now the last claim follows.
\end{proof}

\begin{rema}
Note how the polynomial $h$ determines the geometry of the curve $C$. In fact, the line
$\Hh \cap \Hh'$ is in the support of $C$ if and only if $h = 0$. Furthermore,  $C$ is a
union of two disjoint planar curves if $0 \neq h \in K$. If $h \notin K$, then $C$ is the
union of two planar curves that meet in a zero-dimensional scheme whose degree is less
than each of the degrees of the planar curves.
\end{rema}

\begin{corollary}\label{cor-ModRao-red} Adopt the notation of Lemma \ref{lem-true-curve} and
let  $C \subset \{x y = 0 \}$ be a  curve that is not arithmetically Cohen-Macaulay. Then
the Hartshorne-Rao module of $C$ is:
$$
M_C := H^1_* (\cI_C) \cong \left ({R}/{(x,y,F,G)}\right ) (- d_h).
$$
It is self-dual and $M_C^{\vee} \cong M_C (d-2)$.
\end{corollary}
\begin{proof}
This follows from the Isomorphism \eqref{eq-ext} because $H^1_* (\cI_C) \cong \Ext^2_R
(I_C, R)^{\vee} (4).$
\end{proof}

The first theorem of the introduction follows now easily.

\begin{proof}[Proof of Theorem \ref{thm-red-quad}]
Lemmas \ref{lem-red-quad-necc} and \ref{lem-true-curve} imply that the homogeneous ideal
of $C$ has the required form where $x, y, F, G$ is a regular sequence and $h \in K[z,t]$,
$A \in K[x,z,t]$, $B \in K[y,z,t]$, $F, G \in K[z, t] \setminus K$, and $A B \neq 0$ in
case $h =0$. Given the specific description of the minimal generators, this is equivalent
to the conditions given in the statement.
\end{proof}

%
\section{Curves on a smooth quadric}\label{sec-smoothquadric}

In this section we describe the equations  of curves on a smooth quadric. Without loss of
generality, we consider the quadric $\mathcal{Q}$ defined by $q=xz-yt$.

Since the hyperplane section of $\Qq$ is a divisor in the class $(1, 1)$, any curve on
$\Qq$ is evenly linked to a curve in the class $(d, 0)$ or $(0, d)$ (cf., e.g.\
\cite{KMMNP}). Hence, each minimal curve $C$ on $\Qq$ is in a linear system $|d \f|$
where $d = \deg C$ and $\f$ is a line on $\Qq$. We use this information to explicitly
determine the Hartshorne-Rao module of any curve on $\Qq$. Finally, we
determine the defining equations of the curves in $|d \f|$.

We begin by computing the minimal free resolution of a particular curve in $|d \f|$.

\begin{lem}\label{res_smoothquadric} Let $d \geq 2$ be an integer.
The minimal free resolution of the ideal $I_d = (xz - yt, (x, y)^d)$ is:
\begin{equation}\label{res_ideal}
\small{ \xymatrix{ 0 \hspace{-1.9cm}& \ar[r] &
R^{d-1}(-d-2) \ar[r]_-{{
\begin{bmatrix}
  N_d \\
 -M_d \\
\end{bmatrix}}}
&
 R^{2d}(-d-1)
 \ar[rr]_-{
\begin{bmatrix}
0 & -P_{d-1} \\
  M_{d+1} & N_{d+1} \\
\end{bmatrix}}
&& R(-2) \oplus R^{d+1}(-d) \hspace{-1cm}&
 \ar[r]_{
 \begin{bmatrix}
 q & P_d
 \end{bmatrix}}
& \hspace{-2.9cm} I_{d} \hspace{-3.9cm} & \ar[r] & 0}}
\end{equation}
where, for any $i\geq 2$, $M_i$ and $N_i$ are the  matrices
\begin{equation*}
    M_i=\begin{footnotesize}
    \underbrace{ \begin{bmatrix}
      y &  &   &  0 \\
      -x & y  &   & \vspace{-2mm}  \\
        & -x & \ddots  &   \\
        &   & \ddots  &  y \\
         0 &   &   &   -x\\
    \end{bmatrix}}_{i-1}
    \left.
\begin{array}{c}   \\   \\   \\   \\\\ \end{array}
    \hspace{-5mm} \right\}\tiny{i}\end{footnotesize}
    \hspace{1mm} \text{,} \hspace{5mm}
 N_i=\begin{footnotesize}
 \underbrace{
 \begin{bmatrix}
      z &  &   &  0 \\
      -t & z  &   & \vspace{-2mm}   \\
        & -t & \ddots &   \\
        &   & \ddots &  z \\
        0 &   &    &   -t\\
    \end{bmatrix}}_{i-1}
    \left.
\begin{array}{c}   \\   \\   \\   \\   \\\end{array}
    \hspace{-5mm} \right\}\tiny{i}
    \end{footnotesize} \hspace{1mm} ,
    \end{equation*}
\text{and}
$P_i$ is the $1 \times (i+1)$ matrix
$P_i = \begin{bmatrix}
x^i & x^{i-1} y & \ldots & y^i
\end{bmatrix}$.
\end{lem}
\begin{proof}
Using the relation $N_{i+1}M_i=M_{i+1}N_i$, it is easy to check that
Sequence \eqref{res_ideal} is a complex. Its exactness follows by the
Buchsbaum-Eisenbud criterion.

Alternatively, one can get the free resolution by applying the mapping
cone procedure to the exact sequence that is induced by multiplication
by the quadric $q = xz - yt$:
$$
0\ra R/(x,y)^{d-1}(-2)\stackrel{q}{\ra} R/(x,y)^d
\ra R/I_{d}\ra 0.
$$
\end{proof}

As a consequence, we determine the Hartshorne-Rao modules of curves in
$|d \f|$.

\begin{corollary} \label{cor-hr-mod}
Adopt the assumption and notation of Lemma \ref{res_smoothquadric}. Then the
Hartshorne-Rao module $M_C$ of the curve $C$ defined by the ideal $J =
(xz - yt, (y,
z)^d)$  is:
$$
M_C \cong \coker \left( R^{2d} (-2) \stackrel{\begin{bmatrix}
N_d^t & M_d^t
\end{bmatrix}}{\longrightarrow}  R^{d-1}\right ).
$$
Moreover, its minimal free resolution is of the form
\begin{eqnarray*}
 0 \to R^{d-1} (-d-2) \to R^{2d}(-d-1) \to &  R^{d+1} (-d)      \\
& \oplus  \\
& R^{d+1}(-2) & \to R^{2d} (-1) \to R^{d-1}  \to M_C \to 0.
\end{eqnarray*}
\end{corollary}

\begin{proof}
Let $D$ be the curve defined  by the ideal $I = (xz - yt, (x, y)^d)$. Then the $K$-dual
of its Hartshorne-Rao module is $M_D^{\vee} \cong \Ext^2_R (I, R)(-4)$. Using Lemma
\ref{res_smoothquadric}, we conclude that the minimal free resolution of $M_D^{\vee}(4)$
is of the form:
\begin{eqnarray} \label{dual-res}
\hspace*{0.9cm} 0 \to F_4 \to F_3 \to F_2 \oplus R^{d+1}(d) \to R^{2d} (d+1)
\stackrel{\begin{bmatrix} N_d^t & M_d^t
\end{bmatrix}}{\longrightarrow}  R^{d-1} (d+2) \to M_D^{\vee}(4) \to 0.
\end{eqnarray}
Obviously, there is an inclusion $(xz-yt, y^d) \subset I \cap J$. Since both sides have
the same degree, we get equality, in other words $I$ is geometrically linked to $J$ by
the complete intersection $(xz-yt, y^d)$. Hence it follows (cf., e.g.,
\cite{N-gorliaison}) that $M_C^{\vee} \cong M_D (d-2)$. But the ideal $I$ is transformed
into $J$ by exchanging the variables $x$ and $z$. Thus, the graded Betti numbers of $M_D$
and $M_C$ agree. Dualizing the Resolution \eqref{dual-res}, our claims follow.
\end{proof}

\begin{rema}
(i) Since $M_C$ has exactly $d-1$ minimal generators if $C \in |d
  \f|$, we recover the
fact that each curve in $|d \f|$ is minimal in its even liaison class.

(ii) Using that every curve on $\Qq$ is evenly linked to a curve in
$|d \f|$ for some $d
\geq 1$, we conclude that Corollary \ref{cor-hr-mod} gives a complete
description of the
module structure of curves on a smooth quadric. A first attempt to
achieve such a
classification has been made in \cite{GM-92}, but the results there
are far less
explicit.
\end{rema}

Theorem \ref{thm-dec-smooth} will follow from our next result. Our
original proof was based on Corollary \ref{cor-hr-mod} and complicated.
Discussions with Silvio Greco lead to the much
simpler proof given below.

\begin{lem} \label{lem-key}
Let $C\subset \Qq$ be a curve in $|d\f|$. Then, there are distinct
lines $L_1,\ldots,L_m$
in the same ruling of $\Qq$ as $\f$ such that $C$ has a primary decomposition
$$
I_C= I_{C_1} \cap I_{C_2} \cap \ldots \cap I_{C_m},
$$
where the curve $C_i$ is is defined by
$$
I_{C_i}=I_\Qq+I_{L_i}^{d_i},
$$
with $d_i\geqslant 1$ and $d_1 + \ldots + d_m = d$.
\end{lem}

\begin{proof}
Let $C_1,\ldots,C_m$ be the components of $C$. Since $\f$ is one of
the two free generators of the Picard group of $\Qq$, each curve $C_i$
is linearly equivalent to $d_i \f$ for some $d_i \geq 1$. Moreover,
since by assumption $C_i$ is irreducible, it must be supported on a
line $L_i$. If $d_i = 1$, then $C_i = L_i$. Assume $d_i \geq 2$. Then
it is well-known (cf.\ \cite{GM-92})  that the ideal of $C_i$ is
minimally generated by the quadric defining $\Qq$  and $d_i + 1$
polynomials of degree
$d_i$. It follows that $I_{C_i} \subset I_{\Qq} +
I_{L_i}^{d_i}$. Since both ideals have the same Hilbert fucntion, they
must be equal.
\end{proof}

We conclude by rewriting the ideal of $C$ such that it is generated by
the maximal  minors of a homogeneous matrix.  We write $I_s(A)$ for the ideal
generated by the $s$-minors of matrix $A$. The  result covers Theorem
\ref{thm-dec-smooth}.

\begin{corollary}\label{cor-ideal-smooth}
Let $C\subset \Qq$ be a curve in $|d\f|$. Then
$$
I_C=I_{\Qq} + I_{L_1}^{d_1}  \ldots  I_{L_m}^{d_m}=I_\Qq + I_d(A),
$$
where $A$ is the block diagonal matrix
$$
A=\underbrace{\begin{bmatrix}
  A_1 &  & & \\
   & A_2 & & \\
   &  & \ddots & \\
   & & & A_m \\
\end{bmatrix}}_{d+m}\left.
\begin{array}{c}   \\     \\   \\\\ \end{array}
    \hspace{-5mm} \right\}d
$$
with
$$A_i=
\begin{bmatrix}
  \ell_i & \ell_i' &  &  &  \\
   & \ell_i & \ell_i' &  &  \\
   &  &\ddots & \ddots &  \\
   &  &  & \ell_i & \ell_i' \\
\end{bmatrix} = \begin{bmatrix}
  \,\ell_i\mathcal I_{d_i} & 0 \\
\end{bmatrix}+\begin{bmatrix}
  \,0 &  \ell_i'\mathcal I_{d_i}
\end{bmatrix} \in R^{d_i, d_i},
$$
$d_1 + \ldots + d_m = d$,  and $m$ distinct lines $L_1,\ldots,L_m$ in the same ruling as
$\f$ defined by  $I_{L_i} = (\ell_i, \ell_i')$.
\end{corollary}

\begin{proof}
It suffices to note that
$$
\bigcap_{i=1}^m (I_{\Qq} + I_{L_i}^{d_i}) = I_{\Qq} + \prod_{i=1}^m  I_{L_i}^{d_i}.
$$
\end{proof}

\section{The Hartshorne-Rao modules} \label{sec-HR}

We now describe the possible Hartshorne-Rao modules of curves and the minimal curves on
any quadric. If $C \subset \Qq$ is arithmetically Cohen-Macaulay, then it is in the even
liaison class of a line. Otherwise, we have:

\begin{theorem} \label{thm-hr-sum} Let $D \subset \PP^3$ be a
curve that is not arithmetically Cohen-Macaulay.  Then the minimal curves in the even
liaison class of $D$ lie on a quadric $\Qq$ if and only if the Harts\-horne-Rao module
$M_D$ of $D$ is (up to changes of coordinates and degree shift) presented by one of the
following homogeneous $s \times (2s + 2)$ matrices ($s \geq 1$):
\begin{itemize}
\item[(i)] \begin{equation*}\label{submatrix M} M:= \left[\begin{array}{ccc} x {\mathcal I}_s &
A&\begin{array}{c}
f_1\\ \vdots \\
f_{s}
\end{array}\\
\end{array}
\right],
\end{equation*}
where $I_s (M)$ has codimension 4, $A = (a_{i, j}) \in K[y,z,t]^{s, s+1}$, $\deg f_1 \geq
\deg a_{1,1} -1 + \sum_{j = 1}^s \deg a_{j, j+1}$, and $I_s (A)$ has codimension 2 and is
locally a complete intersection.

\item[(ii)] \begin{equation*} M := \begin{bmatrix}
      y &  -x &   &  & 0 & z & -t &&& 0\\
       & y  & -x  &&&& z & -t \\
        & & & \ddots  &&&&& \ddots   \\
        &  &  \ddots  &&&&& \ddots  \\
         0 &   &   &   y & -x & 0 &&& z & -t\\
    \end{bmatrix}.
\end{equation*}
\end{itemize}

Moreover, the quadric $\Qq$ is uniquely determined by the Hartshorne-Rao module of $D$ if
it has at least three minimal generators, i.e.\ $s \geq 3$. In fact,
if $s \geq 3$ then, in
case {\rm (i)}, the quadric $\Qq$ must be a double plane $2 H$ where
$H$ is defined by
the unique linear form in the annihilator of $M_C$. In case {\rm
  (ii)}, the quadric $\Qq$
must be smooth and is defined by the unique quadratic form in the annihilator of $M_D$.
\end{theorem}

\begin{proof} Since $D$ is not arithmetically Cohen-Macaulay, the quadric $\Qq$ cannot
have rank 3. Thus, the result about the structure of the Hartshorne-Rao module $M_D$
follows by the description of the Hartshorne-Rao modules of curves on a quadric of rank 1
(\cite[Corollary 1.2]{CGN_doubleplane}), of rank 2 (Corollary \ref{cor-ModRao-red}), and
of rank 4 (Corollary \ref{cor-hr-mod}). Note that the Hartshorne-Rao modules of curves on
a reducible quadric are a subclass of those of the curves on a double plane.

Assume now that the Hartshorne-Rao module of $D$ has at least two minimal generators. Let
$C \subset \Qq$ be a minimal curve in the even liaison class of $D$. Then, in case (i),
\cite[Lemma 4.8]{CGN_doubleplane} implies that the annihilator $\Ann_R (M_D)$ contains a
unique linear form $l$ whereas \cite[Theorem 1.1]{CGN_doubleplane} shows that $l^2$ is
the unique quadric in $I_C$ because (with the notation given there) $s \geq 2$ provides
$\deg p \geq 2$.

In case (ii), the  annihilator of $M_C$ does not contain a linear
form, thus $\Qq$ must be smooth.
Furthermore, $\Ann_R (M_C)$ contains a unique quadratic form $q$.
Corollary \ref{cor-hr-mod} implies that $\deg C \geq s+1 \geq 4$, thus Lemma
\ref{res_smoothquadric} provides that $\Qq$ must be defined by $q$.
\end{proof}

Minimal curves are particularly interesting because the Lazarsfeld-Rao property says that
each even liaison class can be recovered from any of its minimal curves by applying
simple operations. Using the notation of the above theorem, we conclude by describing a
minimal curve in each even liaison class of curves on a quadric.

\begin{rema} \label{rem-min-curves}
Let $D \subset \PP^3$ be a curve as in Theorem \ref{thm-hr-sum} and
let $C$ be a minimal
curve in the even liaison class of $D$. In general, $C$ is not unique
and we are going to specify a particular such curve in each class.

If the Hartshorne-Rao module of $D$ is of type (ii), then we can find
a minimal curve $C$ on a smooth quadric and
the homogeneous ideal of $C$ is described in Theorem
\ref{thm-dec-smooth}. In particular, $C$ could be
a multiple line or a union of skew lines.

Assume now that $M_C$ is of type (i). If $M_C$ is not cyclic, then the
proof of Theorem
\ref{thm-hr-sum} provides that $C$ cannot be a reduced curve. Theorem 1.1 in
\cite{CGN_doubleplane} implies that we can choose $C$ as the union of
a double structure
on a planar curve defined by the ideal $(x, p)$ and a planar curve of
degree $\deg f_1 -
\left ( \deg a_{1,1} -1 + \sum_{j = 1}^s \deg a_{j, j+1} \right )$,
where $p$ is the
determinant of the matrix obtained from $A$ by deleting its last column.

Consider now the case $s=1$. Let us rewrite  $M_D$ as $(S/(x, p, F,
G)) (j)$ for some
integer $j$. Choose a homogeneous polynomial $h$ of degree $\deg G -
\deg F - \deg p +1$
such that also $x, p, h F, G$ is a regular sequence. Then the curve
$C$ defined by
\begin{eqnarray*}
I_C & = & (x^2, x p, p^2 h, p h F + x G) \\
& = & (x^2, x p, p^2, p h F + x G) \cap (x, h)
\end{eqnarray*}
is a minimal curve in the class of $D$ of the form described above. If
$p$ has degree
one, then $C$ is also contained in a reducible quadric, thus its
equations must be of the
form described in Theorem \ref{thm-red-quad}. Note that the curves
that are contained in
two quadrics are called extremal curves; they are precisely the curves with the
highest-dimensional Hartshorne-Rao module among all curves with fixed
degree and
arithmetic genus (cf.\ \cite{mdp2} and \cite{N-extremalcurve}). The
only reduced extremal curve is a union of two lines whose
Hartshorne-Rao module is isomorphic to the field $K$, thus it is of type (i)
and of type (ii) in
Theorem \ref{thm-hr-sum}.
\end{rema}

The minimal curves on a reducible quadric are well understood.

\begin{corollary}
Let $C \subset \PP^3$ be a minimal curve that lies on a reducible
quadric. Then $C$ is an extremal curve.
\end{corollary}

\begin{proof}
This follows from Corollary \ref{cor-ModRao-red} and \cite{MDP3}.
\end{proof}
\bigskip

\noindent
{\bf Acknowledgments}
\smallskip

The authors are grateful to Silvio Greco for his helpful comments.

%

\end{document}

\bibitem{N-matematiche}
{\sc U.\ Nagel},
\newblock {\em On the cohomology and genus of projective curves},
\newblock {Le Matematiche {\bf LV}} (2002), 339--351.

\bibitem{R1}
{\sc P.~A.\ Rao},
\newblock {\em Liaison among curves in {$\PP^3$}},
\newblock {Invent.\ Math.\ {\bf 50}}  (1978/79), 205--217.